\documentclass[a4paper,11pt]{article}
\usepackage{amsmath,amssymb,latexsym,amsfonts}
\usepackage{amsthm}

\usepackage [latin1]{inputenc}
\usepackage[english]{babel}

\usepackage{lmodern}

\usepackage{caption}
\usepackage{natbib}
\usepackage{xcolor}
\bibpunct{[}{]}{}{n}{}{}
\usepackage{indentfirst}
\usepackage{upref}

\newtheorem{theorem}{Theorem}
\newtheorem{lemma}[theorem]{Lemma}
\newtheorem{proposition}[theorem]{Proposition}
\newtheorem{corollary}[theorem]{Corollary}

\theoremstyle{definition}

\usepackage{euscript}
\usepackage{yhmath}
\usepackage{enumerate}
\DeclareMathOperator{\diag}{diag}

\def\qed{{$\Box$}}

\newcommand{\ke}{\operatorname{Ker}}

\newcommand{\C}{\ensuremath{\mathbb C}}

\newcommand{\me}{\operatorname{e}}

\newcommand{\al}{\ensuremath{\alpha}}

\newcommand{\la}{\ensuremath{\lambda}}

\newcommand{\si}{\ensuremath{\sigma}}

\newcommand{\La}{\ensuremath{\Lambda}}

\newcommand{\Cnn}{\ensuremath{{\mathbb C}^{n \times n}}}

\newcommand{\Cnuno}{\C^{n\times 1}}

\newcommand{\secsim}[2]{\ensuremath{#1_1,
\ldots,#1_{#2}  }}  

\newcommand{\norma}[1]{\lVert #1 \rVert}

\newcommand{\VS}[2]{\sigma_1(#1)\ge\cdots \ge \sigma_{#2}(#1)}

\newcommand{\VP}[3]{|#1_1(#2)|\ge\cdots \ge |#1_{#3}(#2)|}

\newcommand{\combLineal}[3]{\ensuremath{#1_1 #2_1  +\cdots + {#1}_{#3} {#2}_{#3}}}

\usepackage{ifthen}

\newcounter{natural}

  \renewcommand{\ge}{\geqslant}

\renewcommand{\le}{\leqslant}

\DeclareMathOperator{\dist}{dist}

\usepackage[pagebackref=true,citecolor=red,linkcolor=blue,colorlinks=true,backref=page,urlcolor=blue,plainpages=false,pdfhighlight=/I]{hyperref}

\title{Normal matrices\thanks{This work was supported by the Ministry of Education and Science, Project MTM2017-83624-P.}}
\author{Gorka Armentia, Juan-Miguel Gracia and Francisco-Enrique Velasco\thanks{Department of Mathematics,  University of the Basque Country, Faculty of Pharmacy, Paseo de la Universidad,7, 
01006 Vitoria-Gasteiz, Spain. E-mail addresses: \href{mailto:gorka.armentia@ehu.eus}{gorka.armentia@ehu.eus}, \href{mailto:jm.gracia@telefonica.net}{jm.gracia@telefonica.net}, \href{mailto:franciscoenrique.velasco@ehu.eus}{franciscoenrique.velasco@ehu.eus}
}}
\date{1 June 2021}

\begin{document}
\maketitle
\begin{abstract}
Let $A$ be a square complex matrix and $z$ a complex number. The distance, with respect to the spectral norm, from $A$ to the set of matrices which have $z$ as an eigenvalue is less than or equal to the distance from $z$ to the spectrum of $A$. If these two distances are equal for a sufficiently large finite set of numbers $z$ which are not in the spectrum of $A$, then the matrix $A$ is normal.
\end{abstract}

\section{Distance to the spectrum}
Let $A\in\Cnn$ and $z\in\C$. We let $\La(A)$ stand for the spectrum or set of eigenvalues of $A$. The distance from $z$ to $\La(A)$ is defined by
\[
\dist(z,\La(A)):= \min_{\al \in \La(A)}|z-\al|.
\]
The singular values of $A$ will be denoted by $\VS  A n$. Let $A^*$ denote the conjugate transpose of $A$. We recall that a vector $x\in\Cnuno$ is said to be of \textbf{unit length} if $\norma{x}_2=1$.

The following proposition will be needed to prove the main result of this section.

\begin{proposition}[Weyl]\label{pro:27abril2021-10}
Let $M\in\Cnn$ and let $\la_1(M),\ldots,\la_n(M)$ be the eigenvalues of $M$, which are ordered such that
\[
\VP \la M n.
\]
Then,
\[
\si_1(M)\ge |\la_1(M)|, \qquad     |\la_n(M)|\ge \si_n(M).
\]
\end{proposition}

\begin{proof}
Let $x\in \Cnuno$ be a unit vector such that $Mx=\la_1(M)x$. Then,
\[
\norma{Mx}_2=|\la_1(M)|.
\]
Since
\[
\si_1(M)=\max_{\substack{y\in \Cnuno \\ \norma{y}_2=1}} \norma{My}_2,
\]
we have
\[
\si_1(M)\ge |\la_1(M)|.
\]

\bigskip

Let $U=[\secsim un],V=[\secsim vn]\in \Cnn$ be unitary matrices so that
\[
U^*MV=\diag(\si_1(M),\ldots,\si_n(M)).
\]
Let $x\in \Cnuno$ be a unit eigenvector of $M$ associated with the eigenvalue $\la_n(M)$. As $\{\secsim vn\}$ is an orthonormal basis of $\Cnuno$, there are $\secsim \al n$ $\in \C$ such that
\[
x= \combLineal \al v n
\]
and
\begin{equation}\label{eq:27abril2021-10}
|\al_1|^2+\cdots+|\al_n|^2=1.
\end{equation}
Since $Mx=\la_n(M)x$, we deduce that 
\[
M(\combLineal \al v n)=\al_1 \si_1(M) u_1 +\cdots + \al_n \si_n(M) u_n =\la_n(M)x.
\]
Consequently, the square of the norms of these vectors are equal:
\[
|\al_1|^2 \si_1^2(M)+\cdots +|\al_n|^2 \si_n^2(M)=|\la_n(M)|^2.
\]
But, due to the fact that the singular values are ordered non-increasingly, we deduce that
\[
(|\al_1|^2 +\cdots +|\al_n|^2) \si_n^2(M) \le |\al_1|^2 \si_1^2(M)+\cdots +|\al_n|^2 \si_n^2(M),
\]
and from~\eqref{eq:27abril2021-10}, we obtain
\[
\si_n^2(M)\le |\la_n(M)|^2.
\]
Hence, $\si_n(M)\le |\la_n(M)|$. \hfill\qedhere
\end{proof}

The following result is derived from Proposition~\ref{pro:27abril2021-10}.
\begin{theorem}\label{teo:27abril2021-10}
Let $A\in\Cnn$. Then, for each $z\in \C$,
\begin{equation}\label{eq:27abril2021-20}
\si_n(zI_n-A)\le \dist (z,\La(A)).
\end{equation}
\end{theorem}

\begin{proof}
For each $z\in\C$, let the eigenvalues $\la_1(zI_n-A)$, $\ldots$, $\la_n(zI_n-A)$ of the matrix $zI_n-A$ be arranged non-increasingly according to their moduli:
\[
|\la_1(zI_n-A)|\ge\cdots\ge |\la_n(zI_n-A)|.
\]
As a consequence, we infer that 
\begin{equation}
\vert \la_n(zI_n-A)\vert = \min_{\al\in \La(A)}\vert z-\al\vert. \tag*{\qedhere}
\end{equation}
\let\qed\relax
\end{proof}

\section{A criterion on normality}

Let us recall that a matrix $A\in\Cnn$ is said to be \textbf{normal} if $A^*A=AA^*$. To the extent of our knowledge, around ninety different characterizations of the normal matrices are known. See References~\cite{elsner1998normal} and~\cite{grone1987normal}. We give a different proof of a relatively recent characterization, which appeared as Theorem 2 in~\cite{brooks2018resolvent}. We first show two lemmas.


\begin{lemma}\label{lem:20}
Let $A\in \mathbb{C}^{n\times n}$, $\lambda \in \Lambda(A)$, $x\in \mathbb{C}^n, x\neq 0$,  such that $A x = \lambda x$. Let us assume that
\[
\vert \lambda\vert = \sigma_n(A).
\]
Then,
\begin{enumerate}[{\upshape (i)}]

\item $x^*A = \lambda x^*$.

\item $\lambda$ is a semisimple eigenvalue of $A$.

\end{enumerate}
 
\end{lemma}

\medskip

\begin{proof}
Let us assume that the multiplicity of $\sigma_n(A)$ as singular value of $A$ is $k$; that is
\[
\sigma_n(A)  = \cdots = \sigma_{n-k+1}(A) < \sigma_{n-k}(A).
\]
Consider the singular value decomposition of $A$:
\begin{equation}\label{eq:5}
A v_i = \sigma_i(A) u_i, \hspace*{2cm}  A^* u_i = \sigma_i(A) v_i,
\end{equation}
where $i\in \{1,\ldots,n\}$. Thus, $B_1=\{v_1,\ldots,v_n\}$ and $B_2=\{u_1,\ldots,u_n\}$ are two orthonormal bases for $\mathbb{C}^{n\times 1}$. Without loss of generality, it may be assumed that $x$ is a unit vector. If $x$ is written with respect to the basis $B_1$, then
\begin{equation}\label{eq:7}
x = \sum_{j=1}^{n} \alpha_j v_j,
\end{equation}
with 
\[
\sum_{j=1}^{n} \vert \alpha_j\vert^2 = 1.
\]
On the one hand, $Ax = \lambda x$ by hypothesis. On the other hand, from~\eqref{eq:5} and~\eqref{eq:7} we obtain
\[
A x = A\left(\sum_{j=1}^{n} \alpha_j v_j\right) = \sum_{j=1}^{n} \alpha_j \sigma_j(A) u_j.
\]
Therefore, 
\begin{equation}\label{eq:8}
\lambda x = \sum_{j=1}^{n} \alpha_j \sigma_j(A) u_j.
\end{equation}
By computing the square of the norms of both vectors, we see that
\begin{equation}\label{eq:9}
\vert \lambda \vert^2 = \sum_{j=1}^{n} \vert \alpha_j \vert^2 \sigma^2_j(A).
\end{equation}
Since $\sigma_{n-k+1}(A) = \cdots = \sigma_{n}(A) = \vert \lambda \vert$, we will prove that $\alpha_j=0$, for $j\in \{1,\ldots,n-k\}$. Indeed,~\eqref{eq:9} can be written as
\[
\vert \lambda \vert^2 = \sum_{j=1}^{n-k} \vert \alpha_j \vert^2 \sigma^2_j(A) + \sum_{j=n-k+1}^{n} \vert \alpha_j \vert^2 \vert\lambda\vert^2 ,
\]
and from that
\begin{align*}
\left(1-\sum_{j=n-k+1}^{n} \vert \alpha_j \vert^2\right)\vert \lambda \vert^2& = \sum_{j=1}^{n-k} \vert \alpha_j \vert^2 \sigma^2_j(A) \geqslant \sum_{j=1}^{n-k} \vert \alpha_j \vert^2 \sigma^2_{n-k}(A)\\[0.2cm] & = \left(1-\sum_{j=n-k+1}^{n} \vert \alpha_j \vert^2\right) \sigma^2_{n-k}(A).
\end{align*}
If
\[
1-\sum_{j=n-k+1}^{n} \vert \alpha_j \vert^2 \neq 0,
\]
then $\vert \lambda\vert^2 \geqslant \sigma^2_{n-k}(A)$, which is impossible. Hence, 
\[
1-\sum_{j=n-k+1}^{n} \vert \alpha_j \vert^2 = 0 \Longleftrightarrow \sum_{j=1}^{n-k} \vert \alpha_j \vert^2 = 0
\]
and, consequently, $\alpha_j = 0$, for $j\in \{1,\ldots,n-k\}$. Thus,~\eqref{eq:7} and~\eqref{eq:8} can be expressed as
\begin{equation}\label{eq:20}
x = \sum_{j=n-k+1}^{n} \alpha_j v_j, \hspace*{2cm} x = \me^{-\delta i} \sum_{j=n-k+1}^{n} \alpha_j  u_j,
\end{equation}
where $\me^{\delta i} = \lambda/\vert \lambda \vert = \lambda/\sigma_n(A)$. Finally, using~\eqref{eq:20} in the computation of $x^*A$, we prove the thesis of the lemma.
\begin{align*}
x^* A & = \me^{\delta i} \left(\sum_{j=n-k+1}^{n} \overline{\alpha}_j  u^*_j \right) A = \me^{\delta i} \sum_{j=n-k+1}^{n} \overline{\alpha}_j  \sigma_n(A) v^*_j = \me^{\delta i} \sum_{j=n-k+1}^{n} \overline{\alpha}_j  \vert \lambda \vert v^*_j \\[0.2cm] &= \vert \lambda \vert \me^{\delta i} \sum_{j=n-k+1}^{n} \overline{\alpha}_j v^*_j =\lambda x^*.
\end{align*}

We prove (ii) by contradiction. Let us assume, therefore, that there exists $x\neq 0$ such that not only is an eigenvector of $A$ associated with $\lambda$, but is also the first element of a Jordan chain associated with $\lambda$. The chain is assumed to be of length greater than $1$; that is, there exists a vector $y$ such that
\begin{align}\label{eq:40}
A x &= \lambda x,\notag\\[0.25cm]
A y &= \lambda y + x.
\end{align}
Since $\sigma_n(A) = \vert \lambda\vert$, from (i) we have $x^*A = \lambda x^*$. By multiplying~\eqref{eq:40} on the left by $x^*$, we obtain
\[
x^*A y = \lambda x^*y + x^*x \Longleftrightarrow \lambda x^*y = \lambda x^*y + \Vert x \Vert^2 \Longleftrightarrow \Vert x \Vert^2 =0 \Longleftrightarrow x =0;
\]
which is impossible. \hfill \qedhere
\end{proof}

\bigskip


\begin{lemma}\label{lem:40}
Let $A\in \mathbb{C}^{n\times n}$ and $\Lambda(A)=\{\lambda_1,\ldots,\lambda_p\}$, where $p\leqslant n$. If $x\neq0$ satisfies
\[
A x = \lambda_k x, \hspace*{2cm}  x^*A = \lambda_k x^*,
\]
for some $k\in \{1,\ldots, p\}$, then $x$ is orthogonal to each eigenvector associated with $\lambda_i$, where $i\neq k$.
\end{lemma}

\medskip

\begin{proof}
Let $y\neq 0$ be an eigenvector associated with $\lambda_i$, where $i\in\{1,\ldots,p\}\backslash\{k\}$; that is,
\[
A y = \lambda_i y.
\]
Then,
\[
x^*A y = \lambda_i x^*y.
\]
Since, by hypothesis, $x^*A = \lambda_k x^*$, we have
\[
\lambda_k x^* y = \lambda_i x^*y,
\] 
which implies $x^* y=0$. \hfill \qedhere
\end{proof}

\begin{theorem}\label{teo:11mayo2021-10}
Let $A\in \mathbb{C}^{n\times n}$ and $\Lambda(A)=\{\lambda_1,\ldots,\lambda_p\}$, where $p\leqslant n$. Let us assume that $z_1,\ldots,z_p\in \mathbb{C}\backslash \Lambda(A)$ satisfy
\begin{equation}\label{eq:10}
\sigma_n(z_kI_n-A) = \vert z_k- \lambda_k\vert,
\end{equation}
for $k\in\{1,\ldots,p\}$ respectively. Then $A$ is normal.
\end{theorem}

\medskip

\begin{proof}

Let $k$ be any integer of the set $\{1,\ldots,p\}$. Let us denote by $m_k$ the algebraic multiplicity of $\lambda_k$ and by $s_k$ the geometric multiplicity of $\lambda_k$; that is, 
\[
m_k:= \dim \EuScript{R}_{\lambda_k}(A),\qquad s_k:= \dim \ke(\lambda_kI_n-A);
\]
where $\EuScript{R}_{\lambda_k}(A)$ is the root subspace of $A$ associated with $\lambda_k$. See its definition on page 46 of~\cite{gohberg2006invariant}. Let us assume that the unit vectors $q_{k1},\ldots,q_{k,s_k}$ are linearly independent right eigenvectors associated with $\lambda_k$. Hence, for each $j\in\{1,\ldots,s_k\}$, it is satisfied $A q_{kj} = \lambda_k q_{kj}$, which yields $(z_k I_n - A)q_{kj} = (z_k - \lambda_k) q_{kj}$. Since~\eqref{eq:10} is satisfied, by Lemma~\ref{lem:20} (i) for $z_k-\lambda_k$ and $q_{kj}$, we deduce that $q_{kj}^*(z_kI_n - A) = (z_k-\lambda_k) q_{kj}^*$. Moreover, from Lemma~\ref{lem:20} (ii) it may be concluded that $z_k-\lambda_k$ is a semisimple eigenvalue of $z_k I_n - A$ and, consequently, $\lambda_k$ is a semisimple eigenvalue of $A$. It follows that $\ke(\lambda_kI_n-A) = \ke\left[(\lambda_kI_n-A)^2\right]$; which implies 
\[
\EuScript{R}_{\lambda_k}(A) = \ke(\lambda_kI_n-A)
\]
and, consequently, $m_k=s_k$. By Lemma~\ref{lem:40}, each $q_{kj}$, where $j \in\{1,\ldots,s_k\}$, is orthogonal to $q_{\ell h}$, where $\ell\in\{1,\ldots,p\}\backslash\{k\}$ and $h\in\{1,\ldots,s_\ell\}$. Thus, each $q_{kj}$ is orthogonal to any eigenvector associated with an eigenvalue different from $\lambda_k$, which gives
\[
\EuScript{R}_{\lambda_k}(A)\,\bot\,\EuScript{R}_{\lambda_j}(A)
\]
where $j,k\in\{1,\ldots,p\}$ and $j\neq k$. Therefore, these relations of orthogonality enable us to write the decomposition of $\mathbb{C}^n$ into root subspaces as
\[
\mathbb{C}^n = \EuScript{R}_{\lambda_1}(A) \oplus \cdots \oplus \EuScript{R}_{\lambda_p}(A),
\]
where $\oplus$ stands for orthogonal direct sum. To conclude the proof, it suffices to use the method of Gram-Schmidt in order to build an orthonormal basis of eigenvectors  for each root subspace. Hence, we have proved that there exists an orthonormal basis of $\mathbb{C}^n$ of eigenvectors of $A$. Consequently, by Theorem 4, p. 275 of~\cite{gantmacher1:1966}, $A$ is normal. \hfill \qedhere
\end{proof}

\begin{corollary}
A matrix $A\in\Cnn$ is normal if and only if for each complex number $z$,
\begin{equation*}\label{eq:27abril2021-30}
\si_n(zI_n-A)= \dist (z,\La(A))
\end{equation*}
holds.
\end{corollary}


\end{document}